\documentclass[12pt]{proc-l}

\usepackage{amsmath}
\usepackage{amssymb}
\usepackage{amsthm}
\usepackage{a4wide}
\usepackage[english]{babel}

\numberwithin{equation}{section}

\def\keepspace{\ifnum\catcode`\ =10
  \let\next\keepspacebis \else \let\next\relax \fi
  \next}
\def \keepspacebis{\obeyspaces
  \afterassignment\keepspaceaux\let\next= }
{\obeyspaces%
\gdef\keepspaceaux{%
\ifx \next\space\let\next\ignorespaces\fi%
\catcode`\  =10\relax\next}}

\newcommand{\Z}{\mathbb{Z}}
\newcommand{\Q}{\mathbb{Q}}
\newcommand{\R}{\mathbb{R}}

\def\al{\alpha}

\newcommand{\eps}{\varepsilon}

\newcommand{\moins}{\setminus}

\newcommand{\unr}{\{1, \ldots, r\}}
\newcommand{\uns}{\{1, \ldots, s\}}

\newtheorem{Th}{Theorem}
\newtheorem{Lemme}{Lemma}
\newtheorem{Prop}{Proposition}

\newcommand{\Dem}{\noindent{\bf Proof: }}

\newcommand{\tauxia}{\tau(\xi)}
\newcommand{\tauxib}{\tau'(\xi)}
\newcommand{\tauxic}{\tau''(\xi)}
\newcommand{\tauxid}{\tau'''(\xi)}

\newcommand{\tauxiunsurxizeroa}{\tau(\xi_1/\xi_0)}

\newcommand{\caltxia}{{\mathcal T}(\xi)}
\newcommand{\caltxib}{{\mathcal T}'(\xi)}
\newcommand{\caltxic}{{\mathcal T}''(\xi)}
\newcommand{\caltxid}{{\mathcal T}'''(\xi)}

\newcommand{\dtauxia}{\tau_r(\bfxi)}
\newcommand{\dtauxib}{\tau_r'(\bfxi)}
\newcommand{\dtauxic}{\tau_r''(\bfxi)}
\newcommand{\dtauxid}{\tau_r'''(\bfxi)}
\newcommand{\omzbfxi}{\omega_0(\bfxi)}

\newcommand{\dcaltxia}{{\mathcal T}_r(\bfxi)}
\newcommand{\dcaltxib}{{\mathcal T}_r'(\bfxi)}
\newcommand{\dcaltxic}{{\mathcal T}_r''(\bfxi)}
\newcommand{\dcaltxid}{{\mathcal T}_r'''(\bfxi)}

\newcommand{\dtauxiaunxi}{\tau_1(\bfxi)}
\newcommand{\dtauxiaregaleun}{\tau_1(\xi_0,\xi_1)}
\newcommand{\dtauxibregaleun}{\tau_1'(\xi_0,\xi_1)}
\newcommand{\dtauxicregaleun}{\tau_1''(\xi_0,\xi_1)}
\newcommand{\dtauxidregaleun}{\tau_1'''(\xi_0,\xi_1)}
\newcommand{\omzbfxiregaleun}{\omega_0(\xi_0,\xi_1)}
\newcommand{\omzbfxivraimentr}{\omega_0(\xi_0,\ldots,  \xi_r)}
\newcommand{\omzbfxivraiments}{\omega_0(\xi_0,\ldots,  \xi_s)}

\newcommand{\Span}{{\rm Span}}
 \newcommand{\bfxi}{\underline{\xi}}
  \newcommand{\dtau}{\tau}

\title[Irrationality exponent and rational approximations]
{Irrationality exponent and rational approximations with prescribed growth}

\author{St\'ephane Fischler and Tanguy Rivoal}
\date{June 9, 2009}

\begin{document}


 \begin{abstract}
 Let $\xi$ be a real irrational number. We are interested in sequences of
linear forms in 1 and
$\xi$, with integer coefficients, which tend to $0$. 
Does such a sequence exist 
such that the linear forms are small (with given rate of decrease) and
the coefficients have some
given rate of growth?
If these  rates  are essentially geometric, a necessary condition for  such a sequence
to exist is that
the linear forms are not too small,  a condition which can be expressed
precisely using the
irrationality exponent of $\xi$. We prove that this condition is actually
sufficient, even for
arbitrary rates of growth and decrease. We also make some remarks and ask some questions
about
multivariate generalizations connected to Fischler-Zudilin's new proof of
Nesterenko's linear
independence criterion. 
\end{abstract}
 
\maketitle

\section{Introduction}

In 1978, Ap\'ery~\cite{apery} proved  the irrationality of $\zeta(3)$ by constructing  
two explicit sequences of integers $(u_n)_n$ and $(v_n)_n$ such 
that $0\neq u_n\zeta(3)-v_n\to 0$ and $u_n\to+\infty$, both at geometric rates.  He also  
deduced from this an upper bound for the irrationality 
exponent $\mu(\zeta(3))$ of $\zeta(3)$. In general, the irrationality exponent $\mu(\xi)$ of 
an irrational number $\xi$ is defined as the infimum of all real numbers $\mu$ such that 
the inequality 
$$
\left\vert\xi-\frac pq \right\vert > \frac{1}{q^\mu}
$$
holds for all integers $p,q$, with $q $ sufficiently large. It is well-known that $\mu(\xi)\ge 2$ 
for any irrational number $\xi$ and that it equals $2$ for almost all irrational numbers.
The following lemma is often   used  (as in Ap\'ery's proof) to bound $\mu(\xi)$ from above, for 
example for the numbers $\log(2)$ and $\zeta(2)$. (Other lemmas can be used to bound the 
irrationality exponent of numbers of a different nature, like $\exp(1)$.)

\begin{Lemme} \label{lem1}
Let  $\xi \in \R \moins \Q$, and $\al, \beta$ be real numbers such that $ 0 < \al < 1$ 
and  $  \beta > 1$. Assume there exist integer sequences  $ (u_n)_{n \geq 1} $ 
and  $ (v_n)_{n \geq 1} $ 
such that
\begin{equation} \label{eqlem1}
  \lim_{n\to +\infty}  |u_{n}\xi - v_{n}|^{1/n}  = \al \mbox{ and }   
 \limsup_{n\to+\infty}   |u_{n} |^{1/n}  \leq \beta.
\end{equation}
Then we have 
$\mu(\xi) \leq 1 - \frac{\log \beta}{\log \al}$.
\end{Lemme}

The proof of Lemma~\ref{lem1} is not difficult. Many variants of this result exist; a slightly more 
general version of Lemma~\ref{lem1} will be proved in \S~\ref{subsecvariantelemme}. Another
 variant, proved in~\cite{edd} (Proposition~3.1), asserts that  Lemma~\ref{lem1} holds when~\eqref{eqlem1} 
is replaced with
$$  \limsup_{n\to+\infty}  \frac{|u_{n+1}\xi - v_{n+1}|}{|u_{n}\xi - v_{n}|}  \leq \al
\mbox{ and } 
   \limsup_{n\to+\infty}  \frac{u_{n+1}}{u_{n}}  \leq \beta. $$

\bigskip

In this text, we prove that Lemma~\ref{lem1} and these variants  are best possible, by obtaining a very precise converse result:

\begin{Th} \label{th0}
Let  $\xi \in \R \moins \Q$, and $\al, \beta$ be real numbers such that $ 0 < \al < 1$, 
 $  \beta > 1$ and $\mu(\xi) < 1 - \frac{\log \beta}{\log \al}.$ Then  there exist 
integer sequences  $ (u_n)_{n \geq 1} $ and  $ (v_n)_{n \geq 1} $ 
such that
$$  \lim_{n\to +\infty}  \frac{u_n \xi - v_n}{\al^n} = \lim_{n\to +\infty} \frac{u_n}{\beta^n}= 1$$
and, consequently,
$$ 
  \lim_{n\to +\infty}  \frac{|u_{n+1}\xi - v_{n+1}|}{|u_{n}\xi - v_{n}|}  
= \lim_{n\to +\infty}  |u_{n}\xi - v_{n}|^{1/n}  = \al \mbox{ and } 
 \lim_{n\to +\infty}   \frac{u_{n+1}}{u_{n}} =   \lim_{n\to +\infty}  |u_{n} |^{1/n} = \beta.  $$
\end{Th}

Theorem~\ref{th0} answers completely all questions asked in~\cite{edd}, where 
the {\em density exponent} is defined (see \S~\ref{secedd} below).

\bigskip

An essential feature of Lemma~\ref{lem1}, and all its variants, is that the sequences 
$ (u_n)$ and $(u_n \xi-v_n)$ are assumed to have  essentially geometrical behaviour.  
An assumption of this kind  is necessary, since the convergents of the continued fraction 
expansion of $\xi$ (for instance) always make up a sequence of approximants to $\xi$ that 
are far more precise, but  if $\mu(\xi) > 2$ they don't have    geometrical behaviour. 

\bigskip

However, a   geometrical behaviour is not necessary in Theorem~\ref{th0}, as the following 
generalization  shows.  In this statement, and throughout the paper,   we denote by $o(1)$ any 
sequence that tends  to 0 as $n$ tends to infinity.

\begin{Th} \label{th1} 
Let $\xi \in \R \moins \Q$, and let $(Q_n)$ and $(\eps_n)$ be sequences of positive real numbers with 
$$ \lim_{n\to +\infty} Q_n  = +\infty, \hspace{1cm} \lim_{n\to +\infty} \eps_n = 0 
\hspace{1cm} \mbox{ and} \hspace{1cm} 
\eps_n \geq Q_n^{-\frac{1}{\mu-1}+o(1)},$$
where $\mu$ is a real number such that $\mu > \mu(\xi)$.  

Then there exist   integer sequences $  (u_n) $ and  $  (v_n) $ such that 
$$ \lim_{n\to +\infty} \frac{u_n}{Q_n} = \lim_{n\to +\infty} \frac{ u_{n}\xi - v_{n} }{\eps_n} = 1.$$  
\end{Th}

The important point in this theorem is that our only assumption is that $\eps_n$ is not too 
small, namely
$$\limsup_{n\to+\infty} \frac{-\log \eps_n}{\log Q_n} < \frac{1}{\mu(\xi) -1}.$$

Theorem~\ref{th1} answers the questions asked in \S~8 of~\cite{edd}. 

\bigskip

The structure of this text is as follows. In Section~\ref{secdem}, we prove 
Theorem~\ref{th1} (and therefore, as a special case, Theorem~\ref{th0}). Then we 
recall (in \S~\ref{secedd}) the definition~\cite{edd} of the  {\em density exponent}, and 
deduce from Theorem~\ref{th0} that it is always 0 or $\infty$. In \S~\ref{secexp} 
we prove a slight generalization of Lemma~\ref{lem1} which enables us to obtain a general 
statement (containing Lemma~\ref{lem1}, Theorem~\ref{th0} and Theorem~\ref{th1}) consisting 
in the equality of several exponents of Diophantine approximation. 

Finally we partially 
generalize (in \S~\ref{subsecmulti}) this statement to the multivariate setting, where 
we consider simultaneously 
several real numbers $\xi_0$, \ldots, $\xi_r$ instead of just one $\xi$. 
We then explain the connection between Lemma~\ref{lem1} 
and Nesterenko's linear independence criterion~\cite{Nesterenkocritere}, used in particular 
in the proof (\cite{BR},~\cite{RivoalCRAS}) that $\zeta(s)$ is irrational 
for infinitely many odd integers $s \geq 3$. 
This criterion has been recently generalized in~\cite{SFZu} to take 
into account common divisors to the coefficients of the linear 
forms; the proof involves Minkowski's convex body 
theorem. In  \S~\ref{subsecmulti}, we write down this new proof 
in the case of Nesterenko's criterion only (where Minkowski's convex body 
theorem is simply replaced with Dirichlet's pigeonhole principle) in terms of 
exponents of Diophantine approximation. This enables us to make 
a connection with the other results of the present paper, and to ask 
several questions in the multivariate setting.

\section{Proof of   Theorem~\ref{th1}}  \label{secdem}

The proof of Theorem~\ref{th1} is based on the following lemma, which is  proved 
inside the proof of Lemma~7.3 of~\cite{edd} (p.~39) and is the main step in the 
proof~\cite{edd} that almost all $\xi$ (with respect to Lebesgue measure) have density exponent zero.

 \begin{Lemme} \label{LemmeIMRN}
Let  $c$, $c'$, $\eps$, $Q$  be real numbers such that 
  $1 < c < c' < 2$, $0<\eps < 1$, and  $Q > 1$.

Let  $\xi  $ be an irrational number with $0 < \xi < 1$. Then (at least) one of the following assertions holds:
\begin{enumerate}
\item[$(i)$] There exist coprime integers $u \geq 1$ and $ v \in \{0, \ldots, u\}$ such that
$$u < \frac{2 c^2}{(c-1)(c'-c)} \frac{1}{\eps}$$
and
$$
\Big| \xi -  \frac{v}{u} \Big|  \leq 
\frac{2}{c-1} \left(1 + \frac{c^2}{c'-c}\right) \frac{1}{uQ} .
$$
\item[$(ii)$] There exist integers  $p$  and $q$ such that 
 $$Q \leq q \leq cQ \mbox{ and }   \frac{ \eps}{q} \leq \xi - \frac{p}{q} \leq \frac{ c' \eps}{q}.$$
\end{enumerate}
\end{Lemme}

This lemma is interesting when $\eps $ is much bigger than $1/Q$. It means that, 
unless $\xi$ is very close to a rational number with denominator essentially bounded
 by $1/\eps$, it is possible to find a fraction $p/q$ (which may not be in its lowest 
terms) such that $q$ has essentially the size of $Q$, and $q \xi - p$ that of $\eps$. 
The interesting part, in proving Theorem~\ref{th1}, is that we obtain $ Q \leq q \leq cQ$ 
and $\eps \leq q \xi - p \leq c' \eps$ where $c$ and $c'$ are constants that can be chosen 
arbitrarily close to 1. A variant of this lemma, in which one obtains only 
$ Q \leq q \leq 2Q$ and $\eps \leq q \xi - p \leq 3\eps$, is proved in~\cite{SFrestricted} 
(Lemma 5). The proof uses the same ideas as the one  of 
Lemma~\ref{LemmeIMRN}, but is fairly less complicated.

\bigskip

The proof~\cite{edd} of Lemma~\ref{LemmeIMRN} makes use of Farey fractions. It might be 
possible to prove this lemma using continued fractions.

\bigskip

Let us deduce Theorem~\ref{th1} now.

\begin{proof}
We may assume  $0 < \xi < 1$. Let $(\eta_n)$ be a sequence of positive real numbers 
such that $\lim_{n\to +\infty} \eta_n = 0$ and $\eta_n  = \eps_n^{o(1)}$. We let $\lambda_n = 1+\eta_n$, 
$\mu_n = 1+2\eta_n$, 
$Q'_n = \frac{Q_n}{\sqrt{\lambda_n}}$ and $\eps'_n = \frac{\eps_n}{\sqrt{\mu_n}}$. For 
$n$ sufficiently large, Lemma~\ref{LemmeIMRN} applies with $c =\lambda_n $, 
$c' = \mu_n$, $\eps = \eps'_n $, and $Q = Q'_n $. If $(i)$ holds in this lemma and 
$n$ is sufficiently large,  then we obtain integers $u_n$ and $v_n$ such that 
$$\Big| \xi - \frac{v_n}{u_n}\Big| \leq \frac{2}{\eta_n} 
\Big( 1  + \frac{\lambda_n^2}{\eta_n} \Big) \frac{1}{u_nQ'_n}
\leq  \frac{20}{\eta_n^2 u_n Q_n}  
$$
and
$$u_n <   \frac{2\lambda_n ^2}{\eta_n^2} \frac{1}{\eps'_n} \leq \frac{16 }{\eta_n^2 \eps_n}.
$$ 
Since we have $\eta_n  = \eps_n^{o(1)}$ and  $\eps_n \geq Q_n^{- \frac{1}{\mu-1}+o(1)}$, 
these inequalities yield
$$ \Big| \xi - \frac{v_n}{u_n}\Big| \leq  \frac{20}{\eta_n^2 u_n Q_n}  
\leq \frac{1}{u_n} \Big( \frac{\eps_n \eta_n^2 }{16} \Big) ^{\mu-1+o(1)} 
\leq \frac{1}{u_n^{\mu +o(1)}}$$
which is possible  only for finitely many values of  $n$ since  $\mu > \mu(\xi)$.
Therefore, as soon as $n$ is sufficiently large, Assertion $(ii)$ of 
Lemma~\ref{LemmeIMRN} holds and provides integers $p_n$ and $q_n$ such that 
$$ \frac{Q_n}{\sqrt{\lambda_n}} \leq q_n \leq Q_n \sqrt{\lambda_n} 
\mbox{ and } 
 \frac{\eps_n}{\sqrt{\mu_n}}  \leq q_n \xi - p_n \leq \eps_n \sqrt{\mu_n}.$$
 This concludes the proof of Theorem~\ref{th1}.
\end{proof}

\section{Consequences for the density exponent} \label{secedd}

Let $\xi \in \R \moins \Q$.
For any non-decreasing sequence  ${\bf u}=(u_n)_n$ of positive integers, let us define
$$
\al_\xi({\bf u}):=\limsup_{n} \frac{\vert u_{n+1}\xi-v_{n+1}\vert}
{\vert u_n\xi-v_n\vert},
\quad
\beta({\bf u}):=\limsup_{n} \frac{ u_{n+1}}{u_n},
$$
where $v_n$ is the nearest integer to $u_n\xi$.
We defined in~\cite{edd} the density exponent $\nu(\xi)$ of  $\xi$ as the infimum   of the quantity 
$\log \sqrt{\al_\xi({\bf u})\beta({\bf u})}$ when ${\bf u}$ ranges through the non-decreasing sequences 
such that $\al_\xi({\bf u}) < 1$ and $\beta({\bf u}) <+\infty$ (with the convention 
$\nu(\xi) = +\infty$ if there is no such ${\bf u}$).

We proved in~\cite{edd} that $\nu(\xi)=+\infty$ when $\xi$ is a Liouville number, i.e., 
when  $\mu(\xi)=+\infty$ (that is, when for any $\mu>0$,
there exists a rational number $p/q$ such that $\vert \xi -p/q\vert <1/q^\mu$). 
Theorem~\ref{th1} implies the converse statement, in a more precise form: 

\begin{Th} If  $\xi \in \R \moins \Q$ is not a Liouville number, then $\nu(\xi) = 0$.
\end{Th}

Indeed, we may choose in Theorem~\ref{th0} values of $\al$ and $\beta$ 
arbitrarily close to $1$,  so that the product $\al \beta$ is also arbitrarily close 
to $1$. In a sense, this annihilates the interest of $\nu(\xi)$, since 
it takes only two values ($0$ and $+\infty$) and distinguishes only Liouville numbers
from the other irrational numbers. However, the ideas of~\cite{edd} are at the base 
of the results presented 
in the present paper.

Let us precise here what we expected in~\cite{edd}. We hoped to define a quantity
that would enable us to distinguish between  periods (in the sense of~\cite{KZ}) and  other numbers. 
In particular, we computed upper bounds for $\nu(\xi)$,  for many examples of $\xi$ which are 
periods (see also~\cite{borisedd}). But  
 we did not really take into account another property of the approximations used 
for this: they 
all satisfy a linear recursion of finite order with polynomial coefficients of a special kind. 
Indeed, in all the examples of~\cite{edd}, the sequences $(u_n)_n$ as well as $(v_n)_n$ are such that 
the power series $\sum_{n\ge0} u_n z^n$ and  $\sum_{n\ge0} v_n z^n$ are 
$G$-functions~(\footnote{A power series $\sum_{n\ge0}a_nz^n\in\mathbb{Q}[[z]]$ is a $G$-function 
when: 1) it satisfies a linear differential equation, 2) it has a finite 
 positive radius of convergence, 3) the least commun multiple of the denominators 
of $a_0, a_1, \ldots, a_n$ is bounded by $C^n$ for some $C>0$.}) satisfying the same 
minimal differential equation. This is a very strong property that is not satisfied (in general) by the 
sequences $(u_n)_n$ and $(v_n)_n$ constructed by means of Lemma~\ref{LemmeIMRN} to prove 
Theorem~\ref{th0}.

\section{Exponents of Diophantine Approximation} \label{secexp}

In this section, we state the results of this paper in terms of exponents 
of Diophantine approximation. This enables us to explain the connection with 
Nesterenko's linear independence criterion~\cite{Nesterenkocritere}, and to ask some questions about multivariate 
generalizations of our results.

\subsection{A generalization of Lemma~\ref{lem1}} \label{subsecvariantelemme}

We start with a generalization of the usual Lemma~\ref{lem1}.
We do not write down the proof of this proposition because it 
is a special case of the upper bound 
$  \dtauxia   \leq  \frac{1}{\omzbfxi}$ proved in Theorem~\ref{Thexpmulti} below 
(see \S~\ref{subsecmulti}). 
To deduce Lemma~\ref{lem1} from this  proposition, one takes $\tau = - 
\frac{\log \alpha}{\log \beta}$ and uses the fact that $    \lim_{n\to +\infty} 
\frac{\log  |u_{n}\xi - v_{n}|}{n}  = \log \al$ implies
$  \lim_{n\to +\infty} \frac{\log |u_{n+1}\xi - v_{n+1}|}{\log  |u_{n}\xi - v_{n}|} = 1$.

\begin{Prop}  \label{criterehabituel} 
Let  $\xi \in \R \moins \Q$ and $\tau >0$. Assume    there exist integer sequences $  (u_n) $ 
and  $  (v_n) $ with $u_n \neq 0$ for any $n$, and such that 
$$u_n \xi - v_n \to 0 , \hspace{0.7cm}
|u_{n+1}\xi - v_{n+1}| =  |u_{n}\xi - v_{n}|^{1+o(1)},  \hspace{0.7cm} \mbox{ and} \hspace{0.7cm} 
|u_{n}\xi - v_{n}| \leq |u_n|^{-\tau + o(1)}.$$
Then we have 
$\mu(\xi) \leq 1 + \frac{1}{\tau}$.
\end{Prop}

\subsection{The univariate case}

Let $\xi$ be an irrational real number. Let us consider the following sets:
\begin{itemize}
\item[$\bullet$] $\caltxia$ is the set of all $\tau >0$ for which there exist integer 
sequences $  (u_n) $ and  $  (v_n) $ with $u_n \neq 0$ for any $n$, and 
$$u_n \xi - v_n \to 0 , \hspace{0.7cm}
|u_{n+1}\xi - v_{n+1}| =  |u_{n}\xi - v_{n}|^{1+o(1)},  \hspace{0.7cm} \mbox{ and} \hspace{0.7cm} 
|u_{n}\xi - v_{n}| \leq |u_n|^{-\tau + o(1)}.$$
\item[$\bullet$] $\caltxib$ is the set of all $\tau >0$ for which there exist integer 
sequences $  (u_n) $ and  $  (v_n) $, and $0 < \alpha < 1 < \beta$, with 
$$|u_n \xi - v_n|^{1/n} \to \alpha , \hspace{1cm}
  \limsup_{n\to+\infty}   |u_{n} |^{1/n}  \leq \beta,  \hspace{1cm} \mbox{ and} \hspace{1cm} 
\tau = \frac{-\log \alpha}{\log \beta}. $$
\item[$\bullet$] $\caltxic$ is the set of all $\tau >0$ such that, for any increasing  
sequence $(Q_n)$ of positive integers  there exist   integer sequences 
$  (u_n) $ and  $  (v_n) $ with  
$$   |u_{n}| \leq Q_n^{1+o(1)}  \hspace{1cm} \mbox{ and} \hspace{1cm} 
|u_{n}\xi - v_{n}| = Q_n^{-\tau + o(1)}.$$
\item[$\bullet$] $\caltxid$ is the set of all $\tau >0$ such that, for any   sequences 
$(Q_n)$ and $(\eps_n)$ of positive real numbers with 
$$ \lim_{n\to +\infty} Q_n  = +\infty, \hspace{1cm} \lim_{n\to +\infty} \eps_n = 0 \hspace{1cm} \mbox{ and} \hspace{1cm} 
\eps_n \geq Q_n^{-\tau+o(1)},$$
 there exist   integer sequences $  (u_n) $ and  $  (v_n) $ with  
$$ \lim_{n\to +\infty} \frac{u_n}{Q_n} = \lim_{n\to +\infty} \frac{ u_{n}\xi - v_{n}  }{\eps_n} = 1.$$  
\end{itemize}
Theorem~\ref{Thexp} below shows that   $\tau \leq 1$ for any $\tau $ in $\caltxia$ 
(resp. $\caltxib$, $\caltxic$, $\caltxid$).

We let 
$$\tauxia = \sup \caltxia,$$
and in the same way
$\tauxib = \sup \caltxib,$
$\tauxic = \sup \caltxic,$
$\tauxid = \sup \caltxid,$
with the convention $\sup \emptyset = 0$, so that each of $\tauxia$, $\tauxib$, 
$\tauxic$, $\tauxid$ belongs to $[0,1]$.

If we have $0 < \tau < \tau'$ and $\tau' \in \caltxia$, then $\tau \in \caltxia$ 
so that $\caltxia$ is $ \emptyset $, $(0, \tauxia]$ or $(0, \tauxia)$. The same 
holds for $\caltxib$, $\caltxic$, $\caltxid$.

Moreover the inclusions $\caltxid \subset \caltxic \subset \caltxib \subset 
\caltxia$ hold trivially, so that we have
\begin{equation} \label{eqtri}
\tauxid  \leq \tauxic \leq  \tauxib  \leq  \tauxia.
\end{equation}

\bigskip

The main result of this section is the following 
chain of equalities, which summarizes Lemma~\ref{lem1}, Theorem~\ref{th0} 
and Theorem~\ref{th1}.

\begin{Th} \label{Thexp} For any $\xi \in \R \moins \Q$ we have 
$$\tauxid  = \tauxic = \tauxib =   \tauxia = \frac{1}{\mu(\xi) - 1} \in [0,1]. $$
In particular the following assertions are equivalent:
$\tauxid  =0$; $\tauxic  =0$; $\tauxib  =0$; $\tauxia  =0$; $\xi$ is a Liouville number.
\end{Th}

As a corollary, we have
$ \tauxid  = \tauxic = \tauxib =   \tauxia = 1$ for almost all $\xi$ with respect to 
Lebesgue measure.

\bigskip

\Dem 
Since $2 \leq \mu(\xi) \leq +\infty$ for any $\xi \in \R \moins \Q$, we have 
$\frac{1}{\mu(\xi) - 1} \in [0,1]$.

For any $\tau \in \caltxia$, Proposition~\ref{criterehabituel} yields 
$\mu(\xi) \leq 1 + \frac{1}{\tau}$, that is $\tau \leq \frac{1}{\mu(\xi) - 1}$. 
If $\caltxia \neq \emptyset$ this gives  $\mu(\xi) < \infty$ and 
$ \tauxia \leq \frac{1}{\mu(\xi) - 1}$;  this upper bound holds trivially 
if $\caltxia = \emptyset$.

By Eq.~\eqref{eqtri}, we just have to prove that 
$  \frac{1}{\mu(\xi) - 1} \leq \tauxid $ to finish the proof of 
Theorem~\ref{Thexp}. This is trivial if $\mu(\xi) = + \infty$. Otherwise, 
for any $\mu > \mu(\xi)$, Theorem~\ref{th1} gives $\frac{1}{\mu-1} \in \caltxid$ 
so that $  \frac{1}{\mu(\xi) - 1} \leq \tauxid $. This concludes the proof of Theorem~\ref{Thexp}.

\subsection{The multivariate case} \label{subsecmulti}

Let $\xi_0, \ldots, \xi_r$ be real numbers, with $r \geq 1$. Throughout this section we assume
$$\dim_{\Q} \Span_{\Q} (\xi_0, \ldots, \xi_r) \geq 2$$
so that non-vanishing linear forms in  $\xi_0, \ldots, \xi_r$ with integer 
coefficients can be arbitrarily small. We consider linear forms 
$L = \ell_0 X_0 + \ldots + \ell_r X_r$ with integer coefficients $\ell_i$, and 
we let $H(L) = \max _{0 \leq i \leq r} |\ell_i|$ and $L(\bfxi) =  
\ell_0 \xi_0 + \ldots + \ell_r \xi_r$, where $\bfxi$ stands for the point 
$ (\xi_0, \ldots, \xi_r)$ in $\R^{r+1}$.

 Let us define the following sets:
\begin{itemize}
\item[$\bullet$] $\dcaltxia$ is the set of all $\dtau >0$ for which there 
exists a sequence $(L_n)$ of linear forms  with $L_n(\bfxi) \neq 0$ for any $n$, and 
$$L_n(\bfxi) \to 0 , \hspace{0.7cm}
|L_{n+1}(\bfxi)| =  |L_n(\bfxi)|^{1+o(1)},  \hspace{0.7cm} \mbox{ and} \hspace{0.7cm} 
|L_n(\bfxi)| \leq H(L_n)^{-\dtau + o(1)}.$$
\item[$\bullet$] $\dcaltxib$ is the set of all $\dtau >0$ for which there 
exists a sequence $(L_n)$ of linear forms, and $0 < \alpha < 1 < \beta$, with 
$$|L_n(\bfxi)|^{1/n} \to \alpha , \hspace{1cm}
  \limsup_{n\to+\infty}  H(L_n)^{1/n}  \leq \beta,  \hspace{1cm} \mbox{ and} \hspace{1cm} 
\dtau = \frac{-\log \alpha}{\log \beta}. $$
\item[$\bullet$] $\dcaltxic$ is the set of all $\dtau >0$ such that, for any increasing 
 sequence $(Q_n)$ of positive integers,   there exists  a sequence $(L_n)$ of linear 
forms  with  
$$ H(L_n) \leq Q_n^{1+o(1)}  \hspace{1cm} \mbox{ and} \hspace{1cm} 
|L_n(\bfxi)| = Q_n^{-\dtau + o(1)}.$$
\item[$\bullet$] $\dcaltxid$ is the set of all $\dtau >0$ such that, for any   
sequences $(Q_n)$ and $(\eps_n)$ of positive real numbers with 
$$ \lim_{n\to +\infty} Q_n  = +\infty, \hspace{1cm} \lim_{n\to +\infty} \eps_n = 
0 \hspace{1cm} \mbox{ and} \hspace{1cm} 
\eps_n \geq Q_n^{-\dtau+o(1)},$$
 there exists a sequence $(L_n)$ of linear forms   with  
$$ \lim_{n\to +\infty} \frac{H(L_n)}{Q_n} = \lim_{n\to +\infty} \frac{ L_n(\bfxi)  }{\eps_n} = 1.$$  
\end{itemize}
Theorem~\ref{Thexpmulti} below shows that   $\dtau \leq s$ for any $\dtau $ in 
$\dcaltxia$, with 
 $s = \dim_{\Q} \Span_{\Q} (\xi_0, \ldots, \xi_r) -1$ (and the same holds for  
$\dcaltxib$, $\dcaltxic$ and  $\dcaltxid$).

We let 
$$\dtauxia = \sup \dcaltxia,$$
and in the same way
$\dtauxib = \sup \dcaltxib,$
$\dtauxic = \sup \dcaltxic,$
$\dtauxid = \sup \dcaltxid,$
with the convention $\sup \emptyset = 0$, so that each of $\dtauxia$, $\dtauxib$, 
$\dtauxic$, $\dtauxid$ belongs to $[0,s]$.

If we have $0 < \dtau < \dtau'$ and $\dtau' \in \dcaltxia$, then $\dtau \in \dcaltxia$ 
so that $\dcaltxia$ is $ \emptyset $, $(0, \dtauxia]$ or $(0, \dtauxia)$. The same 
holds for $\dcaltxib$, $\dcaltxic$, $\dcaltxid$.

Moreover the inclusions $\dcaltxid \subset \dcaltxic \subset \dcaltxib \subset 
\dcaltxia$ hold trivially, so that we have
\begin{equation} \label{eqtrimulti}
\dtauxid  \leq \dtauxic \leq  \dtauxib  \leq  \dtauxia.
\end{equation}

Let $\omzbfxi$ be the supremum of the set of all $\omega >0$ such that there exist 
infinitely many $(r+1)$-tuples $(q_0, \ldots, q_r)\in \Z^{r+1}$ with
\begin{equation} \label{eqomun}
| q_i \xi_j - q_j \xi_i | \leq \max(|q_0|, \ldots, |q_r|)^{-\omega}
\mbox{ for any }  1 \leq i < j \leq r.
\end{equation}
Up to renumbering $\xi_0$, \ldots, $\xi_r$, we may assume $\xi_0 \neq 0$ and in 
this case we can replace~\eqref{eqomun} with
\begin{equation} \label{eqomde}
\Big| \frac{\xi_j}{\xi_0} - \frac{q_j}{q_0} \Big| \leq  |q_0| ^{-\omega-1}
\mbox{ for any }  i \in \unr
\end{equation}
so that $\omzbfxi$ measure the quality of simultaneous approximations to 
$(\xi_1/\xi_0, \ldots, \xi_r/\xi_0)$ by rational numbers with the same denominator.

\bigskip

When $r = 1$ and $\xi_0 \neq 0$, we have $\dtauxiaregaleun =  \tauxiunsurxizeroa $ 
(and the analogous equalities for  $\dtauxibregaleun$, $\dtauxicregaleun$ and 
$\dtauxidregaleun$), and $\omzbfxiregaleun = \mu(\xi_1/\xi_0)-1$. 
This explains why the following result is a partial generalization of Theorem~\ref{Thexp}.

\begin{Th} \label{Thexpmulti} 
Let $\xi_0, \ldots, \xi_r \in \R$,  with $r \geq 1$. Then we   have 
$$\dtauxid  \leq \dtauxic   \leq  \dtauxib   \leq   \dtauxia   
\leq  \frac{1}{\omzbfxi}   \leq  s, $$
under the assumption that  $s = \dim_{\Q} \Span_{\Q} (\xi_0, \ldots, \xi_r) -1$ is positive.
\end{Th}

\bigskip

In the rest of this section, we prove this theorem and make some comments.
The upper bounds $\dtauxid  \leq \dtauxic   \leq  \dtauxib   \leq   \dtauxia$ hold 
trivially, and have been observed before in Eq.~\eqref{eqtrimulti}. 

\bigskip

Let us prove that $  \dtauxia   \leq  \frac{1}{\omzbfxi} $. 
As above, we may assume that $\xi_0 = 1$. If $ \dtauxia = 0 $, this result is trivial. 
Otherwise, let $0 < \dtau <  \dtauxia  $ and $0 < \omega <  \omzbfxi  $  
(where $ \omzbfxi  $ could be $+\infty$).  Let $(L_n)$ be a  sequence  of linear 
forms  such that  $L_n(\bfxi) \neq 0$ for any $n$,  and 
\begin{equation} \label{eq20}
L_n(\bfxi) \to 0 , \hspace{1cm}
|L_{n+1}(\bfxi)| =  |L_n(\bfxi)|^{1+o(1)},  \hspace{1cm} \mbox{ and} \hspace{1cm} 
|L_n(\bfxi)| \leq H(L_n)^{-\dtau + o(1)}.
\end{equation}
There exist integers $q_0$, \ldots, $q_r$, with $|q_0|$ arbitrarily large, such 
that~\eqref{eqomde} holds. Let $n$ be the least positive integer such that 
$|q_0 L_n(\bfxi)| \leq 1/2$. Taking  $|q_0|$ sufficiently large ensures that $n$ 
can be made arbitrarily large. Since $n$ is defined in terms of $q_0$,  any number 
denoted by $o(1)$  depends actually on $q_0$, and can be made arbitrarily small by 
choosing $|q_0|$ sufficiently large.

Since $n$ is the least positive integer such that 
$|q_0 L_n(\bfxi)| \leq 1/2$, the integer $n+1$ does not satisfy this property, that is 
 $1/2 < |q_0| |L_{n+1}(\bfxi)| 
=  |q_0|  |L_n(\bfxi)|^{1+o(1)}$,  so that 
\begin{equation} \label{eq21}
|q_0| =  |L_n(\bfxi)|^{-1+o(1)}.
\end{equation}

Now we have (since $\xi_0=1$)
$$L_n(q_0,\ldots,q_r) = q_0 L_n(\bfxi) + L_n(0, q_1-q_0\xi_1, \ldots, q_r-q_0\xi_r).
$$
In the right handside, the first term has absolute value less than or equal to $1/2$, 
by choice of $n$. If the second term has  absolute value less than  the first one, 
then the integer $L_n(q_0,\ldots,q_r)$ 
has  absolute value less than 1 so that it vanishes, and both terms in  the right handside 
have the same   absolute value, thereby contradicting the assumption.

Hence using also~\eqref{eq21},~\eqref{eq20} and~\eqref{eqomde}, we have:
\begin{eqnarray*}
 |L_n(\bfxi)|^{ o(1)} = |q_0 L_n(\bfxi)| &\leq & | L_n(0, q_1-q_0\xi_1, \ldots, q_r-q_0\xi_r)| \\
&\leq &  r H(L_n) \max_{1\leq i \leq r} |q_i - q_0\xi_i| \\
&\leq &  |L_n(\bfxi)|^{ -\frac{1}{\dtau} + o(1)}  |q_0|^{-\omega}
 =   |L_n(\bfxi)|^{\omega -\frac{1}{\dtau} + o(1)} .
 \end{eqnarray*}
Since $\lim_{n\to +\infty} | L_n(\bfxi)| = 0$ and $n$ can be chosen  arbitrarily large, 
this implies $ \omega  \leq 
\frac{1}{\dtau} $. This concludes the proof that 
$  \dtauxia   \leq  \frac{1}{\omzbfxi} $.

\bigskip

Let us prove that $  \frac{1}{\omzbfxi}   \leq  s $, that is $\omzbfxi \geq 1/s$. 
Renumbering $\xi_0$, \ldots, $\xi_r$ if necessary, we may assume that $\xi_0$, \ldots,
 $\xi_s$ are linearly independent over the rationals; then $\xi_{s+1}$, \ldots, 
$\xi_r$ are linear combinations over $\Q$ of these numbers, and it is easy to 
check that $\omzbfxivraimentr = \omzbfxivraiments$. Now the lower bound
 $ \omzbfxivraiments \geq 1/s$ is a classical consequence of Dirichlet's 
pigeonhole principle. Indeed, for any positive integer $Q$,   consider the $Q^s+1$ points 
$(\{q_0 \xi_1/\xi_0\}, \ldots, \{q_0 \xi_s/\xi_0\}) \in [0,1)^s$, for 
$0 \leq q_0 \leq Q^s$ (here $\{x\}$ denotes the fractional part of a 
real number $x$), and the $Q^s$ cubes ${\mathcal C}_{i_1, \ldots, i_s}$ defined, 
for $0 \leq i_1, \ldots, i_s < Q$, by the inequalities $\frac{i_1}{Q} 
\leq x_1 < \frac{i_1+1}{Q}$, \ldots, $\frac{i_s}{Q} \leq x_s < \frac{i_s+1}{Q}$. 
At least two of these points, given (say) by $q'_0$ and $q''_0$, lie in the 
same cube. Letting $q_0 = |q'_0-q''_0|$ and denoting by $q_j$ the nearest integer
to $q_0 \xi_j/\xi_0$ for $j \in \uns$, we obtain~\eqref{eqomde} with 
$\omega = 1/s$. Since a given tuple $(q_0, \ldots, q_s)$ is obtained in this 
way from only finitely many integers $Q$ (because $\xi_0$, \ldots, $\xi_s$ are
 $\Q$-linearly independent), we obtain infinitely many $(s+1)$-tuples 
satisfying~\eqref{eqomde} so that $\omzbfxivraimentr = \omzbfxivraiments \geq 1/s$. 
 This concludes the proof of Theorem~\ref{Thexpmulti}.

\bigskip

A consequence of Theorem~\ref{Thexpmulti} is the inequality $  \dtauxia   
\leq     s$, which amounts to the following statement, known as Nesterenko's 
linear independence criterion~\cite{Nesterenkocritere}:  

\medskip

{\em Assume  there exists a sequence $(L_n)$ of linear forms  with 
$L_n(\bfxi) \neq 0$ for any $n$ and} 
$$L_n(\bfxi) \to 0 , \hspace{1cm}
|L_{n+1}(\bfxi)| =  |L_n(\bfxi)|^{1+o(1)},  \hspace{0.6cm} 
\mbox{ {\em and } } \hspace{0.6cm} 
|L_n(\bfxi)| \leq H(L_n)^{-\dtau + o(1)}$$
{\em  for some  $\dtau >0$. Then we have } 
$\dim_{\Q} \Span_{\Q} (\xi_0, \ldots, \xi_r)  \geq \dtau + 1.$

\medskip

The above arguments provide a simple proof of this criterion (adapted from~\cite{SFZu}), based on 
Dirichlet's pigeonhole   principle and the upper bound $  \dtauxia   
\leq  \frac{1}{\omzbfxi} $ (which is also proved, essentially in the 
same way, as the first step in Nesterenko's inductive proof~\cite{Nesterenkocritere}). 

\medskip

In the one-dimensional case where 
$\bfxi = (1,\xi)$ with $\xi  \in \R \moins \Q$, the upper bound $ \dtauxiaunxi  \leq  
\frac{1}{\omzbfxi} $ corresponds to Proposition~\ref{criterehabituel},  
while   $  \frac{1}{\omzbfxi}   \leq  1 $ simply means $\mu(\xi) \geq 2$ 
(and one way to prove this fact is to  use   Dirichlet's pigeonhole principle 
like in the multivariate setting).

\medskip

It would be very interesting to investigate further around  Theorem~\ref{Thexpmulti}, 
for instance to know for which   $\bfxi$ equality holds 
(as in the univariate case of Theorem~\ref{Thexp}). Is it the case for 
almost all $\bfxi$ with respect to Lebesgue measure? It is well known 
  that $\frac{1}{\omzbfxi}   = s = r$ for almost all $\bfxi$.

\medskip

Another question worth studying is the connection between  $\dtauxia$,     
$\dtauxib$,  $\dtauxic$,  $\dtauxid$, and the exponent $\omega_k(\bfxi)$
 that measures the distance of $\bfxi = (\xi_0, \ldots, \xi_r)$ to subspaces 
of dimension $k+1$ of $\R^{r+1}$ defined over $\Q$, for $k < s$ 
(see~\cite{SchmidtAnnals},~\cite{omegadun},~\cite{omegadde} and~\cite{Nesterenkocritere}).

\medskip

At last, many other questions may be asked about these exponents, 
for instance how to understand the set of 
values taken by $\dtauxia$ (resp. $\dtauxib$,  $\dtauxic$,  $\dtauxid$) as 
$\bfxi$ varies, especially when 
$\bfxi$ is assumed to be of a special form (for instance $\bfxi = (1, \xi, \ldots,  \xi^r)$ with 
$\xi \in \R \moins \Q$). It would be interesting to study if there is any 
connection with Mahler's and Koksma's classifications of numbers.

\providecommand{\bysame}{\leavevmode ---\ }
\providecommand{\og}{``}
\providecommand{\fg}{''}
\providecommand{\smfandname}{\&}
\providecommand{\smfedsname}{\'eds.}
\providecommand{\smfedname}{\'ed.}
\providecommand{\smfmastersthesisname}{M\'emoire}
\providecommand{\smfphdthesisname}{Th\`ese}

\bigskip

S. Fischler,
Univ Paris-Sud, Laboratoire de Math\'ematiques d'Orsay, 
Orsay Cedex, F-91405; CNRS, Orsay cedex, F-91405.

T. Rivoal, Institut Fourier,
CNRS UMR 5582, Universit{\'e} Grenoble 1,
100 rue des Maths, BP~74,
38402 Saint-Martin d'H{\`e}res cedex,
France.

\end{document}